"Workshop on algebraic structures and moduli spaces",
July 14 - 20, 2003, CRM, Univ. de Montreal.

\documentclass[12pt]{amsart}
\usepackage{latexsym,amssymb,amsmath,amsfonts,amsthm}
\usepackage[all]{xypic}

\def\cal{\mathcal}
\def\Bbb{\mathbb}
\def\frak{\mathfrak}
\def\dim{\hbox{dim}\,}
\def\bC{{\Bbb C}}
\def\cpx{\bC}

\def\bT{{\Bbb T}}
\def\bR{{\Bbb R}}

\def\cM{{\cal M}}

\def\cP{{\cal P}}
\def\p{\cP}
\def\cA{{\cal A}}
\def\cG{{\cal G}}
\def\cE{{\cal E}}

\def\cQ{{\cal Q}}
\def\CP{{\Bbb C\Bbb P}}

\def\ss{\sigma}

\def\rk{{\operatorname{rank}}}
\def\Id{{\operatorname{Id}}}

\def\ind{{\operatorname{Ind}}}
\def\End{{\operatorname{End}}}
\def\det{{\operatorname{det}}}

\newtheorem{theorem}{Theorem}

\newtheorem*{remark}{Remark}
\newtheorem{theodef}[theorem]{Theorem-Definition}

\newtheorem*{conjecture}{Conjecture}

\def\prf{\noindent{\em{Proof.}}\rm\ }
\def\endprf{\ \hfill $\Box$\medskip}

\def\hk{hyper\-k\"ahler\ }

\begin{document}

\parindent=0pt

\title{HYPERK\"AHLER NAHM TRANSFORMS}

\author[]{Claudio Bartocci and Marcos Jardim}

\begin{abstract}
Given two \hk manifolds $M$ and $N$ and a quaternionic instanton on their
product, a \hk Nahm transform can be defined, which maps quaternionic
instantons on $M$ to quaternionic instantons on $N$. This construction
includes the case of Nahm transform for periodic instantons on $\bR^4$,
the Fourier-Mukai transform for instantons on K3 surfaces, as well as the
Nahm transform for ALE instantons.
\end{abstract}

\maketitle

\section{Introduction}

The Nahm transform is an instance of ``geometrical Fourier transform'',
originally introduced by W.~Nahm \cite{Nahm} as an extension of the
ADHM method to construct time-invariant instantons (alias {\em monopoles})
on $\bR^4$. This idea was subsequently applied, by P. Braam \& P. van Baal
\cite{BvB}, to study instantons on flat $4$-tori. Nowadays, it is clearly
understood that these constructions are special examples of a much wider
framework, which provides an unified description of correspondences
between solutions of the anti-self-duality equations that are invariant
under dual subgroups of translations of $\bR^4$ (a detailed exposition can
be found in \cite{J5}).

An alternative approach to the Nahm transform is provided by the so-called
Fourier-Mukai transform in algebraic geometry. In this case one obtains
correspondences between moduli spaces of stable sheaves over abelian or K3
surfaces \cite{BBH1}.

Yet another generalization of the Nahm transform can be defined for instantons
over asymptotically locally Euclidean (ALE) 4-manifolds \cite{BaJ}. ALE spaces
are diffeomorphic to minimal resolutions of $\bC^2/\Gamma$, where $\Gamma$ is
a finite subgroup of $SU(2)$, and are endowed with complete \hk
Riemannian metric.
The moduli spaces of instantons over ALE spaces were thoroughly
studied by Nakajima
\cite{Na} and by Kronheimer \& Nakajima \cite{KN}. Contrary to the
previous examples,
the Nahm transform on ALE spaces may fail to be invertible.

It turns out that all these various species of Nahm transforms can be
adequately described in the setting of \hk geometry. Actually, for \hk
manifolds of any dimension, it is possible to generalize the notion of
instanton by defining the so-called quaternionic instantons \cite{MCS, N}.
For these objects, one can prove a generalized Ward correspondence, which
relates quaternionic instantons on $M$ to some holomorphic
vector bundle on the twistor space of $M$. Now, given two \hk manifolds
$M$ and $N$ and a quaternionic instanton $Q$ on $M\times N$, a \hk Nahm
transform is obtained roughly according to the following recipe: take a
quaternionic instanton on $M$  pull it back to $M\times N$, twist it by
$Q$, and push it down to $N$. The existence of the \hk Ward correspondence
ensures (under mild hypotheses) that the result of these operations is a
quaternionic instanton on $N$.

In this expository paper we intend to review, from the unifying \hk viewpoint,
a number of results on Nahm transforms. Section \ref{ward} is devoted to the
description of the Ward correspondence for quaternionic instantons. In Sections
\ref{hknt} and \ref{ex} we present
the \hk Nahm transform, discussing a few fundamental examples which point to an
interesting new conjecture.


\section{Hyperk\"ahler Ward correspondence} \label{ward}

We shall briefly describe the correspondence between quaternionic instantons
on a \hk manifold and certain holomorphic bundles over its twistor space
(for a more general treatment valid for quaternionic K\"ahler manifolds see
\cite{BBH2,MCS}).

Let $(M, g)$ be a (possibly non-compact) $4n$-dimensional \hk manifold (i.e.
the holonomy group $H$ of the Riemannian metric $g$ is contained in $Sp(n)$).
We denote by $\{I_k\}_{k=1,2,3}$ the triple of complex structures and by
$\{\omega_k\}_{k=1,2,3}$ the triple of K\"ahler structures. We recall that
the complex non-degenerate $2$-form $\Omega = \omega_2 + i \omega_3$ is
holomorphic w.r.t.~$I_1$. The essential tool in the study of \hk geometry is
the twistor space, which encodes the information about the $\CP^1$ family of
complex structures on the \hk manifold. The twistor space $Z$ of $M$ is the
sphere bundle associated with the rank $3$ vector bundle $P\times_H
\frak{sp}(1)$,
where $P$ is the holonomy bundle of $M$; we denote by $p: Z \to M$ and by
$q: Z \to \CP^1$ the  natural projections. On $Z$ there is a complex structure
(different from the product structure), and an anti-holomorphic
involution $\tau$,
preserving the fibers of $p$, the so-called twistor lines. The differential of
$\tau$ acts in the following way:
$$d\tau (v, z) = (\bar{v}, d\imath (z)) \quad \forall (v, z)\in T_u Z
= T_{p(u)} M \oplus T_{q(u)} \CP^1\,,$$
where $\imath : \CP^1 \to \CP^1$ is the antipodal map.

For any point $x\in M$, the endomorphism $\Phi= \sum_{k=1}^3 I_k\otimes I_k$ of
$\Lambda^2 T^\ast_x M$ satisfies the identity $\Phi^2 = 2\Phi +
3\Id$, so it has
eigenvalues $3$ and $-1$. The eigenspace $\frak{E}_x$ corresponding to the
eigenvalue $3$ is isomorphic to $\frak{sp}(n)$ and can be characterized as
$$ \frak{E}_x = \bigcup_{u \in p^{-1}(x)} \Lambda_u^{1,1} T^\ast_x M\,,$$
where $\Lambda_u^{1,1} T^\ast_x M$ is the space of the $2$-forms of
type $(1,1)$
w.r.t.~the complex structure associated to the point $q(u)\in \CP^1$.

A {\em quaternionic instanton} on $M$ is a pair $(E, \nabla)$, where $E$
is a complex smooth vector bundle on $M$ and $\nabla$ a connection whose
curvature at any point $x\in M$ takes values in $\frak{E}_x \otimes \End E_x$.
When $\dim M=4$, this definition corresponds to the usual one, since
anti-self-dual
connections are precisely the $2$-forms which are of type $(1,1)$ w.r.t.~all
complex structures compatible with the \hk metric \cite{DK}.

Let $(E, \nabla)$ be a quaternionic instanton on $M$ such that the connection
$\nabla$ is compatible with a given hermitian metric $h$ on $E$. It is easy to
check that the pull-back connection $p^\ast \nabla$ induces a complex structure
on the pull-back bundle $F=p^\ast E$  endowed with the pull-back
hermitian metric
$p^\ast h$. Moreover, we can define a {\em positive real form} $\ss:
F \to F^\ast$,
i.e.~an antilinear anti-holomorphic bundle isomorphism covering the involution
$\tau: Z\to Z$. We set
$$\ss(v)(w) = p^\ast h (v, \bar\tau(w))\,,$$
where $\bar\tau: F \to F$ is the lifting of $\tau$ to $F$. The
following generalized
Ward correspondence holds \cite{BBH2,N}.

\begin{theorem}
Gauge equivalence classes of rank $k$ hermitian quaternionic
instantons on the \hk
manifold $M$ are in one-to-one correspondence with isomorphism
classes of rank $k$
holomorphic bundles on the twistor space $Z$ of $M$, which are holomorphically
trivial along the fibers of $p: Z\to M$ and carry a positive real form.
\end{theorem}

\begin{remark}\rm
When $M$ is compact, irreducible quaternionic $SU(n)$-instantons
correspond to holomorphic
bundles which are $\mu$-stable of degree zero w.r.t.~any K\"ahler
structure induced on $M$
by the \hk structure \cite{V}. In particular, a line bundle fails to
be a quaternionic
instantons if its first Chern class is not orthogonal to the
cohomology classes of the
K\"ahler forms $\{\omega_k\}_{k=1,2,3}$.
\end{remark}


\section{The \hk Nahm transform} \label{hknt}

Let us now suppose we are given the following set of data: two \hk
manifolds $M$, $N$
and an hermitian quaternionic instanton $(Q, \nabla_{Q})$ on the
product $M\times N$.
Let $(Q_\xi,\nabla_{\vert Q_\xi})$ denote the restriction of
$(Q,\nabla_{Q})$ to
$M\times \{\xi\}$, for any $\xi \in N$. If $(E,\nabla)$ is a vector bundle with
connection on $M$, we can define the family of connections
$$ \nabla_\xi =
\nabla \otimes \mathbf{1}_{Q_\xi} + \mathbf{1}_E \otimes
\nabla_{\vert Q_\xi} $$
on the bundle family $E\otimes Q_\xi$, for any $\xi \in N$.

The twistor space $Z_{M\times N}$ of $M\times N$ is the fiber product
along $\CP^1$ of
the twistor spaces  $Z_M$ and $Z_N$ of $M$ and $N$. So, we have the
commutative diagram:
\begin{equation}\label{diagram}
\xymatrix{
Z_M \ar[d]_{p_M} & \ar[l]_{\rho_M}  Z_{M\times N}\ar[d]^q
\ar[r]^{\rho_N} & Z_N \ar[d]^{p_N}\\
M  & \ar[l]^{\pi_M}  {M\times N} \ar[r]_{\pi_N} & N}
\end{equation}
where the horizontal arrows are holomorphic maps while the vertical
arrows are only
smooth maps. By our assumptions, the pull-back of $(Q, \nabla_{Q})$
to $Z_{M\times N}$
is a holomorphic bundle, which is holomorphically trivial along the
fibers of $p:Z\to M$
and carries a positive real form.

Any \hk manifold carries a standard spin structure; the spinor
bundles $S^\pm$ are trivial and the projectivitazion of $S^+$
coincides with the twistor space. By coupling the connection
$\nabla_\xi$ with the Dirac operators on $S^\pm$, we get
the family of Dirac operators:
$$ D_{\xi}^{+} ~ : ~ L^2_{k+1} (E\otimes Q_\xi \otimes S^+) \to
L^2_k (E\otimes Q_\xi \otimes S^-)  $$
and
$$ D_{\xi}^{-} ~ : ~ L^2_k (E\otimes Q_\xi \otimes S^-)  \to
L^2_{k-1} (E\otimes Q_\xi \otimes S^+)  ~, $$
parametrized by the points of $N$. If $M$ is non-compact, one has
to be more careful, and consider suitable weighted Sobolev space
completions  of the spaces of $C^\infty$ sections in order to
guarantee that the Dirac operators $D_{\xi}^{\pm}$ are Fredholm,
but we do not want to enter in these details here (for an example,
see Section \ref{ALE} below).

\begin{theorem} \label{HKNT}
Let $(E,\nabla)$ be an hermitian quaternionic instanton on $M$ and assume
that the index $\hat{E}=\ind(D^\pm)$ of the family of Dirac operators
$D_{\xi}^{\pm}$ is a smooth vector bundle on $N$. Then $\hat{E}\to N$ admits
an hermitian quaternionic instanton $\hat{\nabla}$.
\end{theorem}
\prf The index bundle $\hat{E}$ is endowed with the index connection
$\hat \nabla$, obtained by projecting the trivial connection on
$L^2_k (E\otimes Q_\xi \otimes S^\pm)$ down to $\ker D_\xi^\pm$. It
is easy to check that the pull-back of $(\ind (D^\pm),\hat\nabla)$
to $Z_N$ is a holomorphic bundle holomorphically trivial along the
fibers of $p: Z\to M$ and carrying a positive real form.
\endprf

Alternatively, when $M$ and $N$ are compact, we can use the machinery of higher
direct images of coherent sheaves to reformulate the index
construction above.
Let us fix a complex structure $I_z$ on $M$; we denote by $M_z$ the resulting
complex manifold. If $(E,\nabla)$ is quaternionic instanton, we
denote by $\cE_z$
its associated sheaf of holomorphic sections w.r.t.~the given complex
structure. We
say that $(E,\nabla)$ satisfies the odd (resp.~even) IT condition if
$$ H^k(M, \cE_z\otimes \cQ_{\xi}) = 0 \quad
\hbox{for $k$ even (resp.~odd) and all $\xi \in N$.} $$
The definition is well-posed, since the cohomology groups
$H^k(M,\cE_z\otimes \cQ_{\xi})$
do not depend on the choice of the complex structure \cite{V}.
Consequently, we shall
drop the subscript $z$ in what follows.

By general results (see \cite{BBH}), if $(E,\nabla)$ satisfies the
odd IT condition, then one has the identification of holomorphic vector bundles
$$ - \ind D^+ = \bigoplus_k R^k \pi_{N\ast}(\pi_M^\ast \cE \otimes \cQ)\,.$$
(when the even IT condition is satisfied, one takes $D^-$ instead). Let
$\tilde\cE$ and $\tilde \cQ$ denote the sheaves of holomorphic sections of
the pull-back bundles $p^\ast_M E$ and $q^\ast Q$, respectively. It can be
proved (being careful, since the vertical maps of the diagram \ref{diagram}
are not holomorphic maps) that the pull-back holomorphic bundle
$p^\ast_N(-\ind D^+)$
coincides with $\bigoplus_k
R^k\rho_{N\ast}(\rho_M^\ast\tilde\cE\otimes \tilde \cQ)$.
Summing up, we get the following result.

\begin{theodef}\label{theodef}
Let $M$, $N$ be two compact \hk manifolds, and let $(Q,\nabla_{Q})$
be an hermitian
quaternionic instanton over the product $M\times N$. Assume that
$(E,\nabla)$ is an
hermitian quaternionic instanton on $M$ satisfying the odd IT condition. Then,
$\bigoplus_k R^k \pi_{N\ast}(\pi_M^\ast \cE \otimes\cQ)$ is a
quaternionic instanton
on $M$, which is called the \hk Nahm transform of $(E,\nabla)$ and is
denoted by
$(\hat E,\hat\nabla)$.
\end{theodef}

\begin{remark}\rm
Given the hermitian quaternionic instanton $(Q,\nabla_{Q})$, the correspondence
between quaternionic instantons on $M$ and quaternionic instantons on
$N$ is not
in general invertible. In particular, the \hk Nahm transform of an irreducible
quaternionic instanton may fail to be irreducible.
\end{remark}


\section{Hyperk\"ahler Nahm transforms for instantons} \label{ex}

Let $M$ be a 4-dimensional \hk manifold (not necessarily compact) and let $N$
be a connected component of the moduli of (irreducible) instantons on N (with
some ``framing'' at infinity if $M$ is non-compact); i.e. each point $\xi\in N$
can be regarded as an irreducible anti-self-dual connection $\nabla_\xi$ on a
fixed vector bundle $E\to M$. Then $N$ is also a \hk manifold.

Recall that one can define a {\em universal bundle with connection} over the
product $M\times N$ in the following way \cite{AS,I}. Let $\cA$ denote the set
of all connections on $E$, and let $\cG$ denote the group of gauge 
transformations
(i.e. bundle automorphisms). Again, if $M$ is non-compact, one must 
use the right
analytical framework to define $\cA$ and $\cG$, see Section \ref{ALE} below.
Moreover, let $G$ denote the structure group of $E$,
so that $E$ can be associated with a principal $G_E$-bundle $P$ over 
$M$ by means
of some representation $\rho:G\to\cpx^n$, where $n={\rm rank}~E$. The
gauge group $\cG$ acts on $E\times\cA$ by $g(p,A)=(g(p),g(A))$; This action has
no fixed points,
and it yields a principal $\cG$-bundle $E\times\cA \rightarrow {\cal Q}$,
where ${\cal Q}=E\times\cA/\cG$.

The structure group $G$ also acts on $E\times\cA$, and since this
action commutes
with the one by $\cG$, $G$ acts on $\cal Q$. Moreover, the $G$-action on
${\cal Q}^\ast=E\times\cA^\ast/\cG$ has no fixed points, where
$\cA^\ast$ denotes
the set of irreducible connections on $F$. We end up with a principal
$G$ bundle
${\cal Q}^\ast\rightarrow M\times(\cA^\ast/\cG)$, and we denote by $\tilde{\p}$
the associated vector bundle ${\cal Q}^\ast\times_\rho\cpx^n$. Since $N$ is a
submanifold of $\cA^\ast/\cG$, we define the {\em Poincar\'e bundle}
$\p\to M\times N$
as the restriction of $\tilde{\p}$.

The  principal $G$-bundle ${\cal Q}^\ast$ has a natural connection
$\tilde{\omega}$,
constructed as follows. The space $E\times\cA^\ast$ has a Riemannian
metric which
is equivariant under $G\times\cG$, so that it descends to a
$G$-equivariant metric
on ${\cal Q}^\ast$. The orthogonal complements to the orbits of $G$ yields the
connection $\tilde{\omega}$. Passing to the associated vector bundle
$\tilde{\p}$
and restricting it to $M\times N$ gives a connection $\omega$ on the
{\em Poincar\'e
bundle} $\p$. The pair $(\p,\omega)$ is universal in the sense that
$(\p,\omega)|_{M\times\{\xi\}}\simeq(E,\nabla_\xi)$.

\begin{theorem} \label{instantons}
If $M$ is 4-dimensional \hk manifold and $N$ is a connected component
of the moduli
of instantons on N, then the universal connection $\omega$ on the
Poincar\'e bundle
$\p$ is an hermitian quaternionic instanton.
\end{theorem}
\prf
For any complex structure on $M$ compatible with its \hk metric, the
curvature of the
universal connection is of type $(1,1)$ with respect to induced
complex structure on
the $M\times N$ \cite{I}.
\endprf

As usual, the Weitzenb\"ock formula can be used to guarantee that
$\ker D_\xi^+$ is
trivial for all $\xi\in N$, so that the index bundle is indeed a
smooth vector bundle.
Hence, in view of Theorem \ref{HKNT}, we may define a \hk Nahm
transform that maps
instantons on $M$ into quaternionic instantons on $N$. In particular,
if $N$ is also
4-dimensional, we get a Nahm transform taking instantons on $M$ into
instantons on $N$.

It is worth noting that if $N$ is non-compact, then one must also do
some extra work to determine whether the Nahm transformed quaternionic
instanton has finite energy.

Furthermore, if $M$ can also be regarded as a moduli space of
(quaternionic) instantons
on $N$, one can define a Nahm transform in the reverse direction,
taking (quaternionic)
instantons on $N$ into instantons on $M$. It is then natural to ask
whether such transforms
are the inverse of one another. Based on the known examples of the
\hk Nahm transform
(see below) we pose the following general conjecture.

\begin{conjecture} 
Let $M$ and $N$ be two connected 4-dimensional \hk manifolds. If $N$ is a
component of the moduli of instantons on $M$ and $M$ is a component of the
moduli of instantons on $N$, then:
\begin{itemize}
\item $M$ is diffeomorphic to $N$;
\item the \hk Nahm transform is invertible.
\end{itemize}\end{conjecture}
Whenever $M$ and $N$ are compact and algebraic w.r.t.~at least one
complex structure, then the second claim of the
previous conjecture is actually true, in view of Theorem-Definition
\ref{theodef} and of \cite[Theorem 1.1]{brid}.

\medskip

Let us now analyze some interesting applications of Theorem
\ref{instantons}, mentioning
a few examples in which the conjecture is expected to be true.


\subsection{Algebraic tori}

The simplest example of \hk Nahm transform is for $M$ being an
algebraic torus of complex dimension
$2k$ and $N$ being its dual, regarded as a moduli space of flat
connections on $M$; $\p$
is the usual Poincar\'e bundle. Then the \hk Nahm transform is just
the usual Fourier-Mukai
transform.

In particular, for $M$ being a flat 4-torus, the \hk Nahm transform
coincides with the
usual Nahm transform for instantons on $T^4$ \cite{BvB,DK}, so that
the conjecture is true in this case ($M=T^4$ and $N=(T^4)^*$).


\subsection{K3 surfaces}

The first example  of a {\em non-flat} \hk Nahm transform was described in
\cite{BBH1}. Let $M$ be an algebraic K3 surface, which meets the
following requirements:
\begin{enumerate}
\item $M$ admits a K\"ahler form $\omega$ whose cohomology class $H$
satisfies $H^2=2$;
\item $M$ admits a holomorphic line bundle $L$ whose Chern class
$\ell=c_1(L)$ is such that
$\ell\cdot H=0$ and $\ell^2=-12$;
\item if $D$ is the divisor of a nodal curve on $M$, one has $D\cdot H > 2$.
\end{enumerate}
We say that $M$ is a reflexive K3 surface.

Now let $N$ be the moduli space of instantons of rank 2 with
determinant line bundle $L$
(so that $c_1=\ell$) and $c_2=-1$ over $M$; it can be shown that $N$
is isomorphic to $M$
as a complex algebraic variety \cite{BBH1}. Moreover, any instanton $E$ on $M$,
whose dual
is not a point of $N$, satisfies the odd IT condition. Thus, since
both $M$ and $N$ are
\hk manifolds, we get the following result \cite{BBH1}.

\begin{theorem}
Let $(E, \nabla)$ be an irreducible instanton on $M$ such
that at least one of the following conditions is satisfied:
$\rk E\neq 2$, $\det E\neq  L^\ast$, $c_2(E)\neq 1$.
Then, the Nahm transform of $(E, \nabla)$ is an irreducible instanton
on $N$ having
the same degree. Moreover, the correspondence is invertible.
\end{theorem}

This proves the conjecture for $M$ being reflexive $K3$ surface and
$N$ being the moduli space of instantons on $M$ described above.

It should be pointed out that, also in this case, the \hk Nahm
transform coincides with Fourier-Mukai transform of coherent
sheaves on $M$.


\subsection{ALE spaces} \label{ALE}

Let $\Gamma$ be a finite subgroup of ${\rm SU}(2)$. Let $M$ be the
minimal resolution of the quotient $\bC^2/\Gamma$; it can be proved
that $M$ carries a \hk metric $g$, which is {\em asymptotically locally
Euclidean} (ALE) in the following sense. Some open neighborhood $V$ of
infinity in $M$ has a finite covering $\tilde V$ diffeomorphic to
$\bR^4\backslash \overline{B(0,R)}$, for some $R>0$, and in the induced
coordinates $x_i$ the metric $g$ is required to satisfy the relation
\begin{equation*}
g_{ij} (x) = \delta_{ij} + a_{ij} ~~ {\rm with} ~~
\vert\partial^p a_{ij}(x)\vert = O(\vert x\vert^{(-4-p)}) ~ , ~ p\geq 0
\end{equation*}


The moduli spaces of instantons over ALE spaces have been studied in
detail by Nakajima
\cite{Na}, and by Kronheimer \& Nakajima \cite{KN}. Let $E \to M$ be
a complex (smooth)
vector bundle of rank $n$ and trivial determinant
(i.e. an ${\rm SU}(n)$ vector bundle). In order to define a suitable notion of
connections that are ``framed at infinity'', we fix a group homomorphism
$\rho:\Gamma \to {\rm SU}(n)$; this homomorphism will be identified with a flat
${\rm SU}(n)$ connection over $S^3/\Gamma$. By taking coordinates  $x_i$ on
$V\subset X$ as before, we can extend the function $r(p)=\vert x(p)\vert$
to a positive $r$ function on all of $M$. Given a connection $A_0$, a
weighted Sobolev norm $\Vert \cdot\Vert_{l,2,\delta}$ on the space
of $k$-forms $\Omega^k(E)$ is defined as follows:
\begin{equation}\label{sobolevnorm}
\Vert \alpha \Vert_{l,2,\delta} = \sum^l_{j=0} \Vert r^{j-(\delta+2)}
\nabla_{A_0}^{(j)} \alpha \Vert_{L^2}
\end{equation}
for an integer $l\geq 0$ and $\delta\in \bR$. We denote the completion of
the space $\Omega^k(E)$ in this norm by
$W^{l,2}_{\delta}(E\otimes \Lambda^k T^\ast X)$.

Now fix $l>2$. We say that a connection $A$ on $E$ is
{\it asymptotic} to $\rho$ in $W^{l,2}_{-2}$ if there is a gauge such that
$A=A_0 + \alpha$, where the restriction of $(E,A_0)$ to
$$ \{t\}\times S^3/\Gamma \subset (R,\infty)\times S^3/\Gamma \simeq V $$
is the flat bundle with connection $\rho$, for all $t>R$ and
$$ \Vert \alpha \Vert_{l,2,-2}  <\infty\,. $$
We denote by $\cA^l_X(\rho)$ the space of such connections.

The space of anti-self-dual connections on $E$ asymptotic to $\rho$ and
having topological charge $k$ is described as follows:
$$ \cA^l_{X,asd}(E,\rho,k) =
\left\{ A\in \cA^l_X(\rho) \ \vert \ A \ \text{is anti-self-dual and} \
\frac{1}{8\pi^2} \int_X \Vert F_A \Vert^2=k \right\} $$
The corresponding moduli space is given by the quotient
$$\cM_M(E,\rho,k) =  \cA^l_{X,asd}(E,\rho,k)/ \cG_0^{l+1} $$
where $\cG_0^{l+1}$ is the gauge group of automorphisms of $E$
converging to the
identity, that is:
$$ \cG_0^{l+1} = \left\{ s \in W^{l+1,2}_{-1}({\rm End}(E)) \ \vert \
\Vert s-\mathbf{1}_E \Vert_{l+1,2,-1} < \infty \right\}. $$

The following fundamental result is due to Nakajima \cite{N}:

\begin{theorem}\label{thm1}
Each non-empty, non-compact 4-dimensional component of the moduli space
$\cM_M(E,\rho,k)$ is a complete \hk ALE space.
\end{theorem}

In other words, every such component of the moduli space $\cM_M(E,\rho,k)$
is diffeomorphic to a minimal resolution of $\bC^2/\hat{\Gamma}$ for some
discrete subgroup $\hat{\Gamma}\subset {\rm SU}(2)$, which might be,
in general,
distinct from $\Gamma$.  However, a complete classification (due to Nakajima)
of all possible subgroups $\hat \Gamma$ can be achieved \cite{BaJ2}. Roughly
speaking, Nakajima has shown that the Dynkin diagram associated with
$\hat \Gamma$
has to be a subgraph of the Dynkin diagram associated with $\Gamma$.
As a consequence of this classification, one gets two important
results (the latter was already proved in \cite{KN} by different methods).

\begin{theorem}
Let $M$ and $N$ be two ALE spaces. If $N$ is
diffeomorphic to a 4-dimensional component of the moduli space of
instantons on $M$ and the same holds
for $M$ w.r.t.~$N$, then $M$ is diffeomorphic to $N$.
\end{theorem}

In other words, the first part of the conjecture is true.

\begin{theorem} \label{thm2}
Let $E$ be a rank 2 complex vector bundle over an ALE space $M$. Let
$\iota:\Gamma\hookrightarrow {\rm
SU}(2)$ be the inclusion map, with $|\Gamma|$ denoting the order of
$\Gamma$. Then $\hat{M} = \cM_M(E,\iota,
\frac{|\Gamma|-1}{|\Gamma|})$ is isomorphic to $M$ as a \hk manifold.
\end{theorem}

Let us now consider the \hk manifolds $M$, $N=\hat M$ and an
instanton  $(F,\nabla)$ over $M$.
The following facts are true:
\begin{enumerate}
\item the universal Atiyah-Singer bundle $(Q,\nabla_Q)$ equipped with the
universal connection is a quaternionic instanton by Theorem \ref{instantons}
and \cite{GN};
\item the index $\ind D^-$ of the family of Dirac operators
parametrized by $\hat M$ (cf. Section \ref{hknt}) is a 
finite-dimensional smooth,
vector bundle on $\hat M$, that we shall denote by $\hat F$ \cite{KN};
\item the transformed instanton $\hat\nabla$ has finite energy \cite{BaJ,BaJ2}.
\end{enumerate}

Therefore, \hk Nahm transform takes instantons on $M$ into instantons
on $\hat M$, and vice versa. Using the equivalent Fourier-Mukai formulation,
one can expect to show that the transform is invertible, thus proving also
the second part of the conjecture \cite{BaJ,BaJ2}.


\subsection{Further perspectives}

Let us now consider the case $M=\bT^2\times\bR^2$; let $(x,y)$ denote flat
coordinates in $\bT^2$, and $(r,\theta)$ denote polar coordinates in $\bR^2$.
As it was shown in \cite{BiJ}, the moduli space ${\cal M}_{(k,\xi,\mu,\alpha)}$
of $SU(2)$ instantons $A$ on $M$ which are asymptotic to
$$ A_0=d+i\left(\begin{array}{cc}a_0 & 0 \\ 0 & - a_0
\end{array}\right) ~~ , ~~ {\rm with} $$
$$ a_0 = \lambda_1 dx + \lambda_2 dy + (\mu_1 \cos\theta - \mu_2 \sin\theta)
\frac{dx}{r} + (\mu_1 \sin\theta + \mu_2\cos\theta) \frac{dy}{r} +
\alpha d\theta
$$
is a smooth \hk manifold of real dimension $8k-4$, where
$$ k = \frac{1}{8\pi^2} \int_M |F_A|^2 $$
is the instanton charge. In particular, setting $k=1$, one also shows that
$N={\cal M}_{(1,\xi,\mu,\alpha)} \simeq \bT^2\times\bR^2$ with the flat
metric, whenever it is non-empty.

Thus there is a \hk Nahm transform taking instantons on $M$ into
instantons on $N$; it would be interesting to determine whether the
Nahm transformed instanton satisfy $|F_A|\sim {\rm O}(r^{-2})$ (which
is equivalent to the asymptotic condition above, see \cite{BiJ}).

Moreover, $M$ can also be regarded as a moduli space of instantons on
$N$, so there is a \hk Nahm transform taking instantons on $N$ into
instantons on $M$, giving some evidence for the conjecture being true
when $M=\bT^2\times\bR^2$.

Finally, it would also be interesting to consider the \hk Nahm transform
for instantons over $S^1\times\bR^3$ (alias {\em calorons}) and
$M=\bT^3\times\bR$.
Even more interesting would be to consider different gravitational instantons
(other than ALE spaces), like asymptotically locally flat (ALF) 4-manifolds,
or the ALG gravitational instantons constructed by Cherkis \&
Kapustin in \cite{CK}.
The moduli spaces of instantons over these spaces, however, are much
less understood.


\bigskip
\baselineskip=3pt
{\obeylines
\noindent Claudio Bartocci
\noindent Dipartimento di Matematica
\noindent Universit\`a degli Studi di Genova
\noindent Via Dodecaneso 35
\noindent 16146 Genova, ITALY
\noindent E-mail: bartocci@dima.unige.it
\vskip 1pc
\noindent Marcos Jardim
\noindent Department of Mathematics and Statistics
\noindent University of Massachusetts at Amherst
\noindent Amherst, MA 01003-9305 USA
\noindent E-mail: jardim@math.umass.edu}


\end{document}